\documentclass{amsart}
\usepackage{amssymb,amsmath,amsfonts,amsthm,amscd,latexsym}
\usepackage{amsmath}
\usepackage{latexsym}
\usepackage{graphics}
\usepackage{epsfig}
\usepackage{psfrag}
\usepackage{here}
\usepackage[all]{xy}

\theoremstyle{plain}
\newtheorem*{sat*}{Theorem}
\newtheorem{lemma}{Lemma}[section]
\newtheorem{theorem}[lemma]{Theorem}

\newtheorem{corollary}[lemma]{Corollary}   

\newtheorem{proposition}[lemma]{Proposition}

\theoremstyle{definition}
\newtheorem{definition}[lemma]{Definition}

\newtheorem*{namedtheorem}{\theoremname}
\newcommand{\theoremname}{testing}

\theoremstyle{remark}

\newcommand{\cal}{\mathcal}

\def\_#1{\overline{#1}}

\begin{document}
\title{Geometry
of the complex of curves and of Teichm\"uller space}
\author{Ursula Hamenst\"adt}
\thanks 
{Partially supported by DFG-SPP 1154 and
DFG-SFB 611.}
\maketitle

\begin{abstract} Using train tracks on a nonexceptional
oriented surface $S$ of finite type
in a systematic way we give a proof that the complex of
curves $\mathcal C (S)$ 
of $S$ is a hyperbolic geodesic metric space. We
also discuss the relation
between the geometry of the complex of curves and the geometry
of Teichm\"uller space.
\end{abstract}

\tableofcontents   

\section{Introduction}

Consider a compact oriented surface $S$ of
genus $g\geq 0$ from which $m\geq 0$ points,
so-called \emph{punctures}, have been deleted.
Such a surface is called \emph{of finite type}.
We assume that $S$ is \emph{non-exceptional}, i.e. that 
$3g-3+m\geq 2$; this rules out a
sphere with at most four punctures and a torus
with at most one puncture. 

In \cite{Ha}, Harvey associates to such a surface
the following simplical complex.

\begin{definition}
The \emph{complex of curves} 
\index{complex of curves}\index{complex of
curves} ${\cal C}(S)$
for the surface $S$ is the simplicial complex whose
vertices are the free homotopy classes of
essential simple closed curves on $S$ and whose
simplices
are spanned by collections of such curves which can be realized
disjointly.
\end{definition}
Here we mean by an \emph{essential} simple closed curve
a simple closed curve
which is not contractible nor homotopic into a puncture. 
Since $3g-3+m$ is the
number of curves in a \emph{pants decomposition} of $S$,
\index{pants decomposition} i.e. a maximal 
collection of disjoint mutually not freely homotopic essential
simple closed curves which decompose
$S$ into $2g-2+m$ open subsurfaces homeomorphic to a thrice punctured
sphere, the dimension
of ${\cal C}(S)$ equals $3g-4+m$. 

In the sequel we restrict our attention to
the one-skeleton of the complex of curves 
which is usually called the 
\emph{curve graph}\index{curve graph};
by abuse of notation, we denote it again by ${\cal
C}(S)$. Since $3g-3+m\geq 2$ by assumption,
${\cal C}(S)$ is a nontrivial graph which moreover 
is connected \cite{Ha}.
However, this graph is locally infinite. Namely,
for every simple closed curve $\alpha$ on $S$ the
surface $S-\alpha$
which we obtain by cutting $S$ open along $\alpha$
contains at least one connected component of Euler characteristic
at most $-2$, and such a component contains infinitely many
distinct free homotopy classes of simple closed curves which
viewed as curves in $S$ are disjoint from $\alpha$.

Providing each edge in ${\cal C}(S)$ with the standard
euclidean metric of diameter 1 equips the 
curve graph
with the structure of a geodesic metric space.
Since ${\cal C}(S)$ is not locally finite,
this metric space $({\cal C}(S),d)$
is not locally compact. Masur
and Minsky \cite{MM1} showed that nevertheless its geometry
can be understood quite explicitly. Namely, ${\cal C}(S)$
is hyperbolic of infinite diameter. Here for some $\delta
>0$ a geodesic metric space is
called \emph{$\delta$-hyperbolic in the 
sense of Gromov}\index{$\delta$-hyperbolic}\index{Gromov
hyperbolic} if it
satisfies the \emph{$\delta$-thin triangle condition}: For every
geodesic triangle with sides $a,b,c$ the side $c$ is contained in
the $\delta$-neighborhood of $a\cup b$. Later Bowditch 
\cite{B} gave a
simplified proof of the result of Masur and Minsky which can also
be used to compute explicit bounds for the hyperbolicity constant
$\delta$.

Since the Euler characteristic of $S$ is negative, the surface
$S$ admits a complete hyperbolic metric of finite volume.
The group of diffeomorphisms of $S$ which are isotopic
to the identity acts on the space of such metrics. 
The quotient space under this action is the
\emph{Teichm\"uller space} 
\index{Teichm\"uller space} ${\cal T}_{g,m}$ for $S$
of all  
\emph{marked} isometry classes of 
complete hyperbolic metrics on $S$ of finite
volume, or, equivalently, the space
of all marked complex structures on $S$ of finite type.
The Teichm\"uller space can be equipped with a natural topology,
and with this topology it is homeomorphic to $\mathbb{R}^{6g-6+2m}$.
The \emph{mapping class group} 
\index{mapping class group} ${\cal M}_{g,m}$
of all isotopy classes of orientation preserving diffeomorphisms
of $S$ acts properly discontinuously as a group of diffeomorphisms
of Teichm\"uller space preserving 
a complete Finsler metric, the so-called
\emph{Teichm\"uller metric}\index{Teichm\"uller metric}.
The quotient orbifold is the \emph{moduli space} ${\rm Mod}(S)$
of $S$ of all isometry classes of complete hyperbolic
metrics of finite volume on $S$
(for all this see \cite{IT}).

The significance of the curve graph for the geometry
of Teichm\"uller space comes
from the obvious fact that the mapping class group acts on 
${\cal C}(S)$ as a group of simplicial isometries. Even more
is true: If $S$ is not a twice punctured torus or
a closed surface of genus $2$, then the 
\emph{extended mapping class group} of isotopy classes
of \emph{all} diffeomorphisms of $S$ coincides precisely with
the group of simplicial isometries of ${\cal C}(S)$;
for a closed surface of genus 2, the group of
simplicial isometries of ${\cal C}(S)$ is
the quotient of the extended mapping class group under
the hyperelliptic involution which acts trivially on
${\cal C}(S)$ (see
\cite{I} for an overview on this and related results). 
Moreover, there is a natural map $\Psi:{\cal T}_{g,m}\to
{\cal C}(S)$ which is 
\emph{coarsely ${\cal M}_{g,m}$-equivariant} and 
\emph{coarsely Lipschitz}
with respect to the Teichm\"uller metric on
${\cal T}_{g,m}$.
By this we mean that there is a number $a>1$ 
such that $d(\Psi (\phi h),\phi(\Psi h))\leq a$
for all $h\in {\cal T}_{g,m}$ and all $\phi\in {\cal M}_{g,m}$
and that moreover 
$d(\Psi h,\Psi h^\prime)\leq a d_T(h,h^\prime)+a$
for all $h,h^\prime\in {\cal T}_{g,m}$ where
$d_T$ denotes the distance function
on ${\cal T}_{g,m}$ induced by the Teichm\"uller metric
(see Section 4).

As a consequence,
the geometry of ${\cal C}(S)$ is related to the large-scale
geometry of the Teichm\"uller space and the mapping class group.
We discuss this relation in Section 4. In Section 3
we give a proof of the hyperbolicity
of the curve graph using \emph{train tracks} and splitting
sequences of train tracks in a consistent way as the
main tool. Section 2 introduces 
train tracks, geodesic laminations and quadratic
differentials and summarizes
some of their properties.

\section{Train tracks and geodesic laminations}

Let $S$ be a nonexceptional surface of finite type and
choose a complete hyperbolic metric on $S$
of finite volume. With respect to this metric,
every essential free homotopy class of loops
can be represented
by a closed geodesic which is unique up to parametrization.
This geodesic is \emph{simple}, i.e. without self-intersection,
if and only if the free homotopy class has a simple
representative (see \cite{Bu}). 
In other words, there is a one-to-one correspondence between 
vertices of the curve graph and simple closed geodesics on $S$.
Moreover, there is a fixed \emph{compact} subset $S_0$ of $S$
containing all simple closed geodesics.

The \emph{Hausdorff distance}\index{Hausdorff distance}
between two closed bounded subsets
$A,B$ of a metric space $X$ is defined to be 
the infimum of all numbers $\epsilon >0$ such that
$A$ is contained in the $\epsilon$-neighborhood of $B$ and
$B$ is contained in the $\epsilon$-neighborhood of $A$. 
This defines indeed a distance and hence a 
topology on the
space of closed bounded subsets of $X$; this topology is called
the \emph{Hausdorff topology}\index{Hausdorff topology}. 
If $X$ is compact then
the space of closed subsets of $X$ is compact as well. In particular,
for the distance on $S$ induced by a complete hyperbolic
metric of finite volume, 
the space of closed subsets of 
the compact set $S_0\subset S$
is compact
with respect to the Hausdorff topology.

The collection of all simple closed geodesics on $S$ 
is \emph{not} a closed set with respect to the Hausdorff
topology, but a point in its closure can be described as follows.

\begin{definition}
A \emph{geodesic lamination}\index{geodesic lamination}
for a complete hyperbolic
structure of finite volume on $S$ is a compact subset of $S$ which
is foliated into simple geodesics. 
\end{definition}

Thus every 
simple closed geodesic is a geodesic lamination which consists of a
single leaf. The 
space of geodesic laminations on $S$ 
equipped with the Hausdorff topology is compact,
and it contains the closure of the set of 
simple closed geodesics as a \emph{proper} subset.
Note that every lamination in this closure is necessarily
connected.

To describe the structure of the space of geodesic
laminations more explicitly we introduce some more
terminology.

\begin{definition}
A geodesic lamination $\lambda$ is called
\emph{minimal}\index{minimal geodesic lamination} if
each of its half-leaves is dense in $\lambda$. 
A geodesic lamination $\lambda$ is 
\emph{maximal}\index{maximal geodesic lamination} if all its
complementary components are ideal triangles or once punctured
monogons. A geodesic lamination is 
called \emph{complete}\index{complete geodesic lamination} if it is
maximal and can be approximated in the Hausdorff topology for
compact subsets of $S$ by simple closed geodesics. 
\end{definition}

As an example, a simple
closed geodesic is a minimal geodesic lamination. A minimal
geodesic lamination with more than one leaf has uncountably many
leaves. Every minimal geodesic lamination can be
approximated in the Hausdorff topology by simple closed
geodesics \cite{CEG}. Moreover, a minimal
geodesic lamination $\lambda$ is a \emph{sublamination} 
\index{sublamination} of a
complete geodesic lamination \cite{H1}, i.e. there is a complete
geodesic lamination which contains $\lambda$ as a closed subset.
In particular, every simple closed geodesic on $S$ is a
sublamination of a complete geodesic lamination. 
\emph{Every} geodesic lamination $\lambda$ is a
disjoint union of finitely many minimal components and a finite
number of isolated leaves. Each of the isolated leaves of
$\lambda$ either is an isolated closed geodesic and hence a
minimal component, or it \emph{spirals} about one or two minimal
components \cite{Bo, CEG, O}. This means that the set
of accumulation points of an isolated half-leaf of $\lambda$
is a minimal component of $\lambda$.

Geodesic laminations which are disjoint unions of 
minimal components 
can be equipped with the following additional structure.

\begin{definition} 
A \emph{measured geodesic lamination}\index{measured
geodesic lamination} 
is a geodesic lamination
together with a translation invariant 
\emph{transverse measure}\index{transverse measure}.
\end{definition}

A transverse measure for a geodesic lamination 
$\lambda$ assigns to every smooth 
compact arc $c$ on $S$ with
endpoints in the complement of $\lambda$ and which
intersects $\lambda$ transversely a finite Borel 
measure on $c$ supported
in $c\cap \lambda$. These measures transform in the
natural way under homotopies of $c$ by smooth arcs
transverse to $\lambda$ which move the endpoints of the
arc $c$ within fixed complementary components
of $\lambda$.
The support of the measure is the smallest sublamination
$\nu$ of $\lambda$ such that the measure on any such arc
$c$ which does not intersect $\nu$ is trivial.
This support is necessarily a union of minimal components of 
$\lambda$. An example for a measured geodesic lamination
is a \emph{weighted simple closed geodesic} which consists
of a simple closed geodesic $\alpha$ and a positive weight
$a>0$. The measure disposed on a transverse
arc $c$ is then the sum of the Dirac masses 
on the intersection points between $c$ and $\alpha$ multiplied
with the weight $a$.

The space ${\cal M\cal L}$ of measured geodesic
laminations on $S$ can naturally be equipped with the 
\emph{weak$^*$-topology}. This to\-po\-lo\-gy re\-stricts to
the weak$^*$-to\-po\-lo\-gy on the space of measures on 
a given arc $c$ which is transverse to each lamination from
an open subset 
of lamination space. 
The natural action of the group $(0,\infty)$ by
scaling is continuous with respect to this topology, and
the quotient is the space ${\cal P\cal M\cal L}$ of 
\emph{projective measured geodesic laminations}. This space is
homeomorphic to a sphere of dimension $6g-7+2m$
(see \cite{CEG, FLP, PH}).

The \emph{intersection number} $i(\gamma,\delta)$ between two
simple closed curves $\gamma,\delta\in {\cal C}(S)$ equals the
minimal number of intersection points between representatives of
the free homotopy classes of $\gamma,\delta$. This intersection
function extends to a continuous pairing 
$i:{\cal M\cal L}\times{\cal M\cal L}\to [0,\infty)$,
called the 
\emph{intersection form}\index{intersection form}.

Measured geodesic laminations are intimately related to
more classical objects associated to Riemann surfaces,
namely \emph{holomorphic quadratic 
differentials.}\index{holomorphic 
quadratic differential}
A holomorphic quadratic differential $q$ 
on a Riemann surface $S$ assigns to each complex
coordinate $z$ an expression of the form $q(z)dz^2$
where $q(z)$ is a holomorphic function on the domain
of the coordinate system, and $q(z)(dz/dw)^2=q(w)$
for overlapping coordinates $z,w$. We require
that $q$ has at most a simple pole at each puncture of $S$.
If $q$ does not
vanish identically, then its zeros are isolated and
independent of
the choice of a complex coordinate. If $p\in S$ is not a zero
for $q$ then there is a coordinate $z$ near $p$, unique
up to multiplication with $\pm 1$, such that
$p$ corresponds to the origin and that $q(z)\equiv 1$.
Writing $z=x+iy$ for this coordinate, the euclidean metric
$dx^2+dy^2$ is uniquely determined by $q$. The arcs parallel
to the $x$-axis (or $y$-axis, respectively) define
a foliation ${\cal F}_h$ (or ${\cal F}_v$) on the set of regular
points of $q$ called
the \emph{horizontal} (or \emph{vertical})
foliation. The vertical length $\vert dy\vert$
defines a \emph{transverse measure} for the horizontal
foliation, and the horizontal length $\vert dx\vert$
defines a transverse measure for the vertical foliation.
The foliations ${\cal F}_h,{\cal F}_v$ have
singularities of the same type at the zeros of $q$
and at the punctures of $S$
(see \cite{S} for more on quadratic
differentials and measured foliations).

There is a 
one-to-one correspondence between measured geodesic
laminations and (equivalence classes of) measured foliations
on $S$ (see \cite{L} for a precise statement).
The pair of measured foliations defined by 
a quadratic differential $q$ corresponds under this
identification to a pair
of measured geodesic laminations $\lambda\not=\mu\in {\cal M\cal L}$ 
which \emph{jointly fill
up $S$}. This means that for every $\eta\in {\cal M\cal L}$
we have $i(\lambda,\eta)+i(\mu,\eta)>0$.

Vice versa, every pair of measured geodesic laminations
$\lambda\not= \mu\in {\cal M\cal L}$ which jointly fill
up $S$ defines a unique complex structure of finite type on $S$
together with a holomorphic quadratic differential
$q(\lambda,\mu)$ (see \cite{Ke}
and the references given there) whose \emph{area},
i.e. the area of the singular
euclidean metric defined by $q(\lambda,\mu)$, equals
$i(\lambda,\mu)$. If $\alpha,\beta$ are 
\emph{simple multi-curves} on $S$, which means
that $\alpha$ and $\beta$ consist of 
collections $c=c_1\cup \dots \cup c_\ell\subset
{\cal C}(S)$ of free homotopy classes of 
simple closed curves which can be realized
disjointly, and if $\alpha,\beta$ jointly
fill up $S$, then for all
$a>0,b>0$ the quadratic
differential $q(a\alpha,b\beta)$ defined
by the measured geodesic laminations
$a\alpha,b\beta$ can explicitly be constructed
as follows. Choose smooth
representatives of $\alpha,\beta$, again denoted
by $\alpha,\beta$, which intersect transversely in precisely
$i(\alpha,\beta)$ points; for example, the geodesic
representatives of $\alpha,\beta$ with respect to any complete
hyperbolic metric on $S$ of finite volume have this property. 
For each intersection point between $\alpha,\beta$ choose a closed
rectangle in $S$
with piecewise smooth boundary containing this point in
its interior and which does not contain any other intersection
point between $\alpha,\beta$. We allow that some of the
vertices of such a rectangle are punctures of $S$.
These rectangles can be chosen in
such a way that they provide $S$ with the structure of a cubical
complex: The boundary of each component $D$ of $S-\alpha-\beta$
is a polygon with an even number of sides which 
are subarcs of $\alpha,\beta$ in alternating order.
If $D$ does not contain a puncture, then its boundary has
at least four sides.
Thus we can construct the rectangles in such a way that their
union is all of $S$ and that the
intersection between any two distinct such rectangles either is a
common side or a common vertex. Each rectangle from this cubical
complex has two sides which are ``parallel'' to $\alpha$ and two
sides ``parallel'' to $\beta$
(see Section 4 of \cite{Ke} for a detailed discussion of
this construction).

Equip each rectangle with an euclidean
metric such that the sides parallel to $\alpha$ are of length $b$,
the sides parallel to $\beta$ are of length $a$ and such that the
metrics on two rectangles coincide on
a common boundary arc. These metrics define a piecewise
euclidean metric on $S$ with a
singularity of cone angle $k\pi\geq 3\pi$
in the interior of each disc component
of $S-\alpha-\beta$ whose boundary 
consists of $2k\geq 6$ sides.
The metric also has a singularity
of cone angle $\pi$ at 
each puncture of $S$ which is contained in
a punctured disc component with two sides. 
Since
there are precisely $i(\alpha,\beta)$ rectangles, the area of
this singular euclidean
metric on $S$ equals $ab i(\alpha,\beta)$. 
The line segments of this singular
euclidean metric which are parallel to $\alpha$ and $\beta$ define
singular foliations ${\cal F}_{\alpha},{\cal F}_{\beta}$ on 
$S$ with transverse measures induced by the singular
metric. The metric defines a complex structure on $S$
and a quadratic differential
$q(a\alpha,b\beta)$ which is holomorphic for this structure
and whose horizontal and vertical foliations are just
${\cal F}_\alpha,{\cal F}_{\beta}$, 
with transverse measures determined
by the weights $a$ and $b$. The assignment which
associates to $a>0$ the Riemann surface structure determined
by $q(a\alpha,\beta)$ 
is up to parametrization the
geodesic  in the Teichm\"uller space with respect to
the \emph{Teichm\"uller metric} whose cotangent bundle
contains the differentials $q(a\alpha,\beta)$
(compare \cite{IT, Ke}).

A \emph{geodesic segment} for a quadratic 
differential $q$ 
is a path in $S$ not containing any singularities in its
interior and which is a geodesic in the local euclidean 
structure defined by $q$. A closed geodesic is composed of
a finite number of such geodesic segments which meet at
singular points of $q$ and make an angle at least $\pi$
on either side. Every essential closed curve
$c$ on $S$ is freely homotopic to a closed geodesic with respect
to $q$, and the length of such a geodesic $\eta$ 
is the infimum
of the $q$-lengths of any curve freely homotopic to $\eta$
(compare \cite{R} for a detailed discussion of the
technical difficulties caused by the punctures of $S$)
and will be called the \emph{$q$-length} of our
closed curve $c$. If $q=q(\lambda,\nu)$ for
$\lambda,\nu\in {\cal M\cal L}$ then 
this $q$-length is bounded from
above by $2i(\lambda,c)+2i(\mu,c)$ (see \cite{R}).

Thurston invented a 
way to understand the struc\-ture of the space of
geo\-de\-sic la\-mi\-na\-tions by squeezing almost parallel
strands of such a lamination to a simple arc and 
analyzing the resulting graph. The structure of such a
graph is as follows.

\begin{definition}
A \emph{train track}\index{train track}
on the surface $S$ is an embedded
1-complex $\tau\subset S$ whose edges
(called \emph{branches}) are smooth arcs with
well-defined tangent vectors at the endpoints. At any vertex
(called a \emph{switch}) the incident edges are mutually tangent.
Through each switch there is a path of class $C^1$
which is embedded
in $\tau$ and contains the switch in its interior. In
particular, the branches which are incident
on a fixed switch are divided into
``incoming'' and ``outgoing'' branches according to their inward
pointing tangent at the switch. Each closed curve component of
$\tau$ has a unique bivalent switch, and all other switches are at
least trivalent.
The complementary regions of the
train track have negative Euler characteristic, which means
that they are different from discs with $0,1$ or
$2$ cusps at the boundary and different from
annuli and once-punctured discs
with no cusps at the boundary.
A train track is called \emph{generic} if all switches are at most
trivalent.
\end{definition}

In the sequel 
we only consider generic train tracks.
For such a train track $\tau$, every complementary component
is a bordered subsurface of $S$ whose boundary
consists of a finite number of arcs of class $C^1$ which
come together at a finite number of \emph{cusps}. 
Moreover, for every switch of $\tau$ there is precisely
one complementary component containing the switch in
its closure which has a cusp at the switch.
A detailed
account on train tracks can be found in \cite{PH} and 
\cite{M}.

A geodesic lamination or a train track $\lambda$ is \emph{carried}
by a train track $\tau$ if there is a map
$F:S\to S$ of class $C^1$ which is isotopic to the identity and
maps $\lambda$ to $\tau$ in such a way that the restriction of its
differential $dF$ to every tangent line of $\lambda$ is
non-singular. Note that this makes sense since a train track has
a tangent line everywhere. 

If $c$ is a simple closed curve carried by $\tau$ with carrying
map $F:c\to \tau$ then $c$ defines a \emph{counting measure}
$\mu_c$ on $\tau$. This counting measure is the non-negative
weight function on the branches of $\tau$ which associates
to an open branch $b$ of $\tau$ the number of connected components
of $F^{-1}(b)$. A counting measure is an example for a
\emph{transverse measure}\index{transverse measure} 
on $\tau$ which is defined to be a
nonnegative weight function $\mu$ on the branches of $\tau$
satisfying the \emph{switch condition}: for every switch $s$ of
$\tau$, the sum of the weights over all incoming branches at $s$
is required to coincide with the sum of the weights over all
outgoing branches at $s$. The set $V(\tau)$ of all transverse
measures on $\tau$ is a closed convex cone in a linear space and
hence topologically it is a closed cell. 
More generally, every
measured geodesic lamination $\lambda$ on $S$ which is carried
by $\tau$ via a carrying map $F:\lambda\to \tau$ defines a transverse
measure on $\tau$ by assigning to a branch $b$ the total mass of 
the pre-image of $b$ under $F$; the resulting weight function is
independent of the particular choice of $F$.
Moreover, every transverse measure for $\tau$ can be obtained
in this way (see \cite{PH}).

\begin{definition}
A train track is called
\emph{recurrent}\index{recurrent train track}
if it admits a transverse measure which is
positive on every branch. 
A train track $\tau$ is called \emph{transversely
recurrent}\index{transversely recurrent train track}
if every branch $b$ of $\tau$ is intersected by an
embedded simple closed curve $c=c(b)\subset S$ which intersects
$\tau$ transversely and is such that $S-\tau-c$ does not contain
an embedded \emph{bigon}, i.e. a disc with two corners at the
boundary. 
A recurrent and transversely recurrent train track
is called \emph{birecurrent}. 
A generic transversely recurrent
train track which carries a complete geodesic
lamination is called \emph{complete}\index{complete train track}.
\end{definition}

For every recurrent train
track $\tau$, measures which are positive on every branch
define the interior of the convex
cone $V(\tau)$ of all transverse measures. 
A complete train track is birecurrent \cite{H1}.

A half-branch $\tilde b$ in a generic
train track $\tau$ incident on a switch $v$
is called
\emph{large} if the switch $v$ is trivalent and if
every arc $\rho:(-\epsilon,\epsilon)\to
\tau$ of class $C^1$ which passes through $v$ meets the
interior of $\tilde b$.
A branch $b$ in $\tau$ is called
\emph{large} if each of its two half-branches is
large; in this case $b$ is necessarily incident on two distinct
switches
(for all this, see \cite{PH}).

There is a simple way to modify a transversely
recurrent train track $\tau$ to
another transversely recurrent 
train track. Namely, if $e$ is a large branch of
$\tau$ then we can perform a right or left
\emph{split} of $\tau$
at $e$ as shown in Figure A below.
The split $\tau^\prime$ of a
train track $\tau$ is carried by $\tau$.
If $\tau$ is complete and if the complete geodesic
lamination $\lambda$ is carried by $\tau$,
then for every large
branch $e$ of $\tau$ there is a unique choice of a right or left 
split of $\tau$ at $e$ with the property that the split track
$\tau^\prime$ carries $\lambda$, and $\tau^\prime$ is complete.
In particular, a complete train track $\tau$ can always
be split at any large branch $e$ to a complete train track
$\tau^\prime$; however there may be a choice of a right
or left split at $e$ such that the resulting train track
is not complete any more (compare p.120 in \cite{PH}).

\begin{figure}[ht]
\centering
\includegraphics{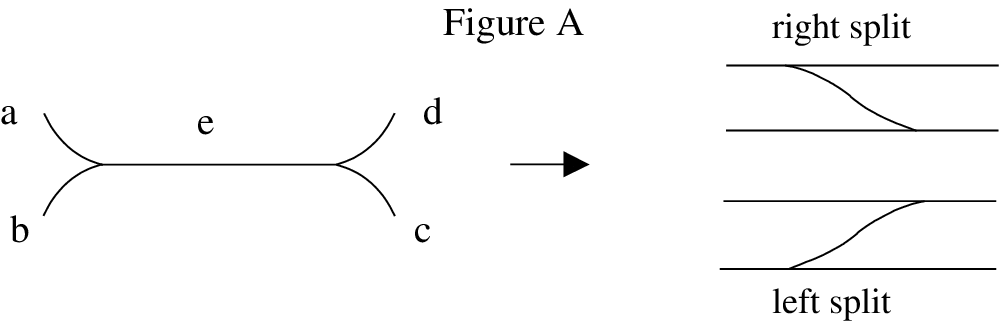}
\end{figure}

In the sequel we denote by ${\cal T}T$ the collection of all
isotopy classes of complete train tracks on $S$. A sequence
$(\tau_i)\subset {\cal T}T$ of complete train tracks is called a
\emph{splitting sequence}\index{splitting sequence} 
if $\tau_{i+1}$ can be obtained from
$\tau_i$ by a single split at some large branch $e$.

\section{Hyperbolicity of the complex of curves}

In this section we present a proof of hyperbolicity of the
curve graph using the main strategy of Masur 
and Minsky \cite{MM1} and 
Bowditch \cite{B} in a modified form.
The first step consists in
guessing a family of uniform 
\emph{quasi-geodesics}\index{quasi-geodesic} in the
curve graph connecting any two points.
Here a $p$-quasi-geodesic for some $p>1$
is a curve $c:[a,b]\to {\cal C}(S)$ which satisfies
\[d(c(s),c(t))/p-p\leq \vert s-t\vert \leq p d(c(s),c(t))+p
\quad \hbox{for all}\, s,t\in [a,b].\]
Note that a quasi-geodesic does not have to be continuous.
In a hyperbolic geodesic metric space, every 
$p$-quasi-geodesic
is contained in a fixed tubular neighborhood of any geodesic
joining the same endpoints, so
the $\delta$-thin triangle condition also holds for
triangles whose sides are
uniform quasi-geodesics \cite{BH}.
As a consequence, for every
triangle in a hyperbolic geodesic 
metric space with uniform quasi-geodesic sides
there is a ``midpoint'' whose distance
to each side of the triangle is bounded from above
by a universal constant.
The second step of the proof consists in finding such a midpoint
for triangles whose sides are curves
of the distinguished curve family.
This is then
used in a third step
to establish the $\delta$-thin triangle condition
for the distinguished family of curves and derive from
this hyperbolicity of ${\cal C}(S)$.
By abuse of notation, in the sequel we simply
write $\alpha\in {\cal C}(S)$ if $\alpha$ is a free
homotopy class of an essential simple closed curve on $S$,
i.e. if $\alpha$ is a vertex of ${\cal C}(S)$.

We begin with defining a map from the set
${\cal T}T$ of complete train tracks on $S$ into
${\cal C}(S)$. For this we call
a transverse measure $\mu$ for a complete train track $\tau$
a \emph{vertex cycle}\index{vertex cycle} 
\cite{MM1} if $\mu$ spans an extreme ray
in the convex cone $V(\tau)$ of all transverse measures on $\tau$.

Up to scaling, every vertex cycle $\mu$ is a counting measure of a
simple closed curve $c$ which is carried by $\tau$ \cite{MM1}. 
Namely, the switch conditions are a family of 
linear equations with integer coefficients for the transverse
measures on $\tau$. Thus an extreme ray is spanned by
a nonnegative \emph{rational} solution which can be scaled 
to a nonnegative integral solution. 
From every 
integral transverse measure $\mu$ for $\tau$ we can construct
a unique \emph{simple weighted multi-curve}, i.e.
a simple multi-curve together with a family
of weights for each of its components, which is
carried by $\tau$ and whose
counting measure coincides with $\mu$ as follows.
For each branch $b$ of $\tau$ draw $\mu(b)$ disjoint arcs 
parallel to $b$. By the switch condition, 
the endpoints of these arcs can
be connected near the switches in a unique way so that
the resulting family of arcs does not have self-intersections.
Let $c$ be the simple multi-curve consisting of the free homotopy
classes of
the connected components of the resulting curve $\tilde c$.
To each such homotopy class associate the number of components
of $\tilde c$ in this class as a weight. The resulting
simple weighted multi-curve is carried by $\tau$, and its counting
measure equals $\mu$. Thus if there are at least 
two components of $\tilde c$
which are not freely homotopic then 
the weighted counting measures of these components 
determine a decomposition of 
$\mu$ into transverse measures for $\tau$ 
which are not multiples of $\mu$. This is impossible
if $\mu$ is a vertex cycle. Hence $c$ consists of 
a single component and  
up to scaling, $\mu$ is the counting measure of a 
simple closed curve on $S$.

A simple closed curve which is carried by $\tau$, with carrying map
$F:c\to \tau$, defines a vertex cycle for $\tau$ only if
$F(c)$ passes through every branch of $\tau$ at most twice, with
different orientation (Lemma 2.2 of \cite{H2}). 
In particular, the counting measure $\mu_c$ of 
a simple closed curve $c$ which defines
a vertex cycle for $\tau$ satisfies
$\mu_c(b)\leq 2$ for every branch $b$ of $\tau$.

In the sequel we mean by a vertex cycle of a complete train track
$\tau$ an \emph{integral} transverse measure on $\tau$ which is
the counting measure of a simple closed curve $c$ on $S$ carried
by $\tau$ and which spans an extreme ray of $V(\tau)$; we also use
the notion vertex cycle for the simple closed curve $c$.
Since the number of branches of a complete
train track on $S$ only depends on the
topological type of $S$, the
number of vertex cycles for a complete train track on $S$ is
bounded by a universal constant (see \cite{MM1} and \cite{H2}).

The following observation of Penner and Harer
\cite{PH} is essential for all what follows. Denote by
${\cal M\cal C}(S)$ the space of all simple
multi-curves on $S$. Let  
$P=\cup_{i=1}^{3g-3+m}\gamma_i\in {\cal M\cal C}(S)$
be a pants decomposition
for $S$, i.e. a simple multi-curve with the maximal number
of components.
Then there is a special family of complete train tracks
with the property  that each pants curve $\gamma_i$ admits a
closed neighborhood $A$ diffeomorphic to an annulus and such that
$\tau\cap A$ is diffeomorphic to a \emph{standard twist connector}
depicted in Figure B. 
\begin{figure}[hb]
\centering
\includegraphics{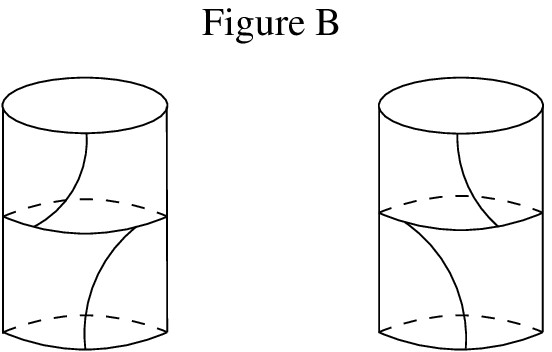}
\end{figure}
Such a train track clearly carries each
pants curve from the pants decomposition $P$ as a vertex cycle; we
call it \emph{adapted} to $P$.
For every complete geodesic lamination
$\lambda$ there is a train track $\tau$ adapted to $P$ which
carries $\lambda$ (\cite{PH}, see also
\cite{H1},\cite{H2}).

Since every simple multi-curve is a subset of a pants
decomposition of $S$, we can conclude.

\begin{lemma}[\cite{H4}]
For every pair
$(\alpha,\beta)\in {\cal M\cal C}(S)\times {\cal M\cal C}(S)$
there is a
splitting sequence $(\tau_i)_{0\leq i\leq m}\subset {\cal T}T$ of
complete train tracks with the property that
$\tau_0$ is adapted to a pants decomposition $P_\alpha\supset \alpha$
and that each component of $\beta$ is a vertex cycle for $\tau_m$.
\end{lemma}

We call a splitting sequence as in the lemma an $\alpha\to
\beta$-splitting sequence. Note that such a sequence is by no
means unique.

The distance in ${\cal C}(S)$ between
two simple closed curves $\alpha,\beta$ is bounded from above by
$i(\alpha,\beta)+1$ 
(Lemma 1.1 of \cite{B} and Lemma 2.1 of \cite{MM1}). In
particular, there is a number $D_0>0$ with the following
property. Let $\tau,\tau^\prime\in {\cal T}T$ and assume
that $\tau^\prime$ is obtained from $\tau$ by at most one
split. Then 
the distance in ${\cal
C}(S)$ between any vertex cycle of $\tau$ 
and any vertex cycle of $\tau^\prime$ is at most $D_0$
(see \cite{MM1} and the discussion following Corollary 2.3 in \cite{H2}).

Define a map $\Phi:{\cal T}T\to {\cal C}(S)$ by assigning to a
train track $\tau\in {\cal T}T$ a vertex cycle $\Phi(\tau)$ for
$\tau$. By our above discussion, for any two choices
$\Phi,\Phi^\prime$ of such a map we have
$d(\Phi(\tau),\Phi^\prime(\tau))
\leq D_0$ for all $\tau\in {\cal T}T$.
Images under the map $\Phi$ of
splitting sequences then define a family of curves in
${\cal C}(S)$ which connect any pair of points in a $D_0$-dense
subset of ${\cal C}(S)\times {\cal C}(S)$, equipped with the
product metric. As a consequence, we can use such
images of splitting
sequences as our guesses for uniform quasi-geodesics.
It turns out that up to parametrization, these curves
are indeed $p$-quasi-geodesics in ${\cal C}(S)$ for a universal
number $p>0$ only depending on the topological type of the
surface $S$ (\cite{MM3}, see also \cite{H2}).

To explain this fact we use
the following construction of Bowditch \cite{B}.
For multi-curves $\alpha,\beta \in
{\cal M\cal C}(S)$ which jointly fill up $S$, i.e. which
cut $S$ into components which are homeomorphic
to discs and once punctured discs, 
and for a number $a>0$ 
let $q(a\alpha,\beta/ai(\alpha,\beta))$ be
the area one quadratic differential whose
horizontal foliation corresponds to the measured geodesic 
lamination $a\alpha$ and whose vertical
measured foliation corresponds to the measured
geodesic lamination $\beta/ai(\alpha,\beta)$.
For $r>0$ define
\[L_a(\alpha,\beta,r)=\{\gamma\in {\cal C}(S)\mid
\max\{ai(\gamma,\alpha),i(\gamma,\beta)/ai(\alpha,\beta)\}\leq r\}.\]
Then $L_a(\alpha,\beta,r)$ is contained
in the set of all all simple closed curves
on $S$ whose $q(a\alpha,\beta/ai(\alpha,\beta))$-length 
does not exceed $2r$. 
Note that we have
$L_a(\alpha,\beta,r)=L_{1/ai(\alpha,\beta)}(\beta,\alpha,r)$
for all $r>0$, moreover
$\alpha^\prime\in L_a(\alpha,\beta,r)$ for every component
$\alpha^\prime$ of $\alpha$ and
every
sufficiently large $a>0$, and $\beta^\prime\in L_a(\alpha,\beta,r)$
for every component $\beta^\prime$ of $\beta$ and
every sufficiently small $a>0$.
Thus for fixed $r>0$ we can
think of a suitably chosen assignment which associates to a
number $s>0$ a point in $L_s(\alpha,\beta,r)$ as a curve in ${\cal
C}(S)$ connecting a component of
$\beta$ to a component of $\alpha$ (provided, of course, that
the sets $L_s(\alpha,\beta,r)$ are non-empty).
Lemma 2.5 of \cite{H2} 
links such curves to splitting sequences.

\begin{lemma}[\cite{H2}]
There is a number $k_0\geq 1$ with the following
property. Let $P$ be a pants
decomposition of $S$, let $\alpha\in {\cal M\cal C}(S)$
be such that $\alpha$ and $P$ jointly fill up $S$ and let 
$(\tau_i)_{0\leq i\leq m}\subset {\cal T}T$ be a
$P\to \alpha$-splitting sequence.
Then there is a non-decreasing
surjective function $\kappa:(0,\infty)\to \{0,\dots,m\}$ such that
$\kappa(s)=0$ for all sufficiently small $s>0, \kappa(s)=m$ for
all sufficiently large $s>0$ and that
for all $s\in (0,\infty)$ there is a vertex cycle
of $\tau_{\kappa(s)}$ which is contained in
$L_s(\alpha,P,k_0)$.
\end{lemma}

Since for every multi-curve $\alpha\in {\cal M\cal C}(S)$
and every pants decomposition $P$ of $S$
there is
a $P\to \alpha$-splitting sequence, we conclude
that for
every $k\geq k_0$ and every $s>0$
the set $L_s(\alpha,P,k)$
is non-empty.
To obtain a control of the size of these sets, 
Bowditch \cite{B} uses the following 
observation (Lemma 4.1 in \cite{B}) whose first part
was earlier shown by 
Masur and Minsky (Lemma 5.1 of \cite{MM1}).

\begin{lemma}[\cite{B}]
There is a number $k_1\geq k_0$
with the following
property. For all $\alpha,\beta\in {\cal M\cal C}(S)$ 
which jointly full up $S$ and every $a\in
(0,\infty)$ there is some $\delta\in L_a(\alpha,\beta,k_1)$ such
that for every $\gamma\in {\cal M\cal C}(S)$ we have
\[i(\delta,\gamma)\leq
k_1\max\{ai(\alpha,\gamma),i(\gamma,\beta)/ai(\alpha,\beta)\}.\] 
In particular, for every $R>0$, for all $\alpha,\beta\in
{\cal M \cal C}(S)$ and for every $a>0$
the diameter of the set $L_a(\alpha,\beta,R)$
is not bigger than $2k_1R+1$.
\end{lemma}

\begin{proof} In \cite{MM1, B} it is shown
that there is a
number $\nu >0$ only depending on the topological type of $S$ and
there is an embedded essential annulus in $S$ whose width
with respect to the piecewise euclidean metric 
defined by the quadratic 
differential $q(a\alpha,\beta/ai(\alpha,\beta))$ 
is at least $\nu$. This means 
that the distance between the boundary circles
of the annulus is at least $\nu$. Assuming the
existence of such an annulus, let $\delta$ be its core-curve.
Then for every simple closed curve $\gamma$ on $S$ and for every
essential intersection of $\gamma$ with $\delta$ there is a subarc
of $\gamma$ which crosses through this annulus and hence whose length
is at least $\nu$; moreover, different subarcs of $\gamma$ corresponding
to different essential intersections between $\gamma$ and $\delta$ are
disjoint. Thus the length with respect to the singular
euclidean metric on $S$ of any simple closed curve $\gamma$ on $S$
is at least $\nu i(\gamma,\delta)$.
On the other hand, by construction the minimal length
with respect to this metric
of a curve in the free
homotopy class of $\gamma$ is bounded from
above by $2\max\{ai(\alpha,\gamma),i(\beta,\gamma)/ai(\alpha,\beta)\}$ and
therefore the core curve $\delta$ of the annulus
has the
properties stated in the first part of our lemma (see \cite{B}).

The second part of the lemma is
immediate from the first. Namely, let $\alpha,\beta\in {\cal
M \cal C}(S)$ and let $a>0$. Choose $\delta\in L_a(\alpha,\beta,k_1)$
which satisfies the properties stated in the first part of the
lemma. If  $\gamma\in L_a(\alpha,\beta,R)$ for some $R>0$
then we have
$i(\gamma,\delta)\leq k_1 R$ and hence $d(\gamma,\delta)\leq
k_1R+1$. \end{proof}

As an immediate consequence of Lemma 3.2 and Lemma 3.3
we observe that there is a universal number $D_1>0$ with the
following property.
Let $P$ be a pants decomposition for $S$, let
$\beta\in
{\cal M\cal C}(S)$ and let $(\tau_i)_{0\leq i\leq m}\subset
{\cal T}T$
be any $P\to \beta$-splitting sequence;
then the Hausdorff distance in ${\cal
C}(S)$ between the sets $\{\Phi(\tau_i)\mid 0\leq i\leq m\}$ and
$\cup_{a>0}L_a(\beta,P,k_1)$ is at most $D_1/16$.
If $c>0$ and if $j\leq m$ is such that there is a vertex
cycle $\gamma$ for $\tau(j)$ which is contained
in $L_c(\beta,P,k_1)$ then the splitting
sequence $(\tau_i)_{0\leq i\leq j}$ is a
$P\to \gamma$-splitting sequence and hence
the Hausdorff distance between $\cup_{a>0}L_a(\gamma,P,k_1)$
and $\cup_{a\geq c}L_a(\beta,P,k_1)$ is at most $D_1/8$.
Moreover, 
for every $\beta\in {\cal C}(S)$ and every simple
multi-curve $Q$ containing $\beta$ as a component the
Hausdorff-distance between the sets
$\cup_{a>0}L_a(\beta,P,k_1)$ and $\cup_{a>0}L_a(Q,P,k_1)$ is at
most $D_1/8$. Thus if $Q,Q^\prime$ are pants
decompositions for $S$ containing a common curve
$\beta\in {\cal C}(S)$ then the Hausdorff distance
between $\cup_aL_a(Q,P,k_1)$ and
$\cup_aL_a(Q^\prime,P,k_1)$ is at most $D_1/4$.

On the other hand, for multi-curves $P,Q\in {\cal M\cal C}(S)$
we have 
\[\cup_{a>0}L_a(P,Q,k_1)=\cup_{a>0}L_a(Q,P,k_1).\]
Therefore from two applications of our above consideration
we obtain the following. Let
$\alpha,\beta\in {\cal C}(S)$ and let
$P,P^\prime,Q,Q^\prime$ be any pants decompositions for $S$
containing $\alpha,\beta$; then
the Hausdorff distance between
$\cup_{a>0}L_a(P,Q,k_1)$ and
$\cup_{a>0}L_a(P^\prime,Q^\prime,k_1)$ is not
bigger than $D_1/2$.
By our choice of $D_1$ this implies that the Hausdorff distance
between the images
under $\Phi$ of \emph{any} $\alpha\to\beta$- or
$\beta\to\alpha$-splitting
sequences is bounded from above by $D_1$.

Now let $\alpha,\beta,\gamma\in {\cal C}(S)$ be such that their
pairwise distance in ${\cal C}(S)$ is at least 3;
then any two of these curves jointly fill up $S$.
Choose pants decompositions $P_\alpha,P_\beta,P_\gamma$ containing
$\alpha,\beta,\gamma$.
Then there are
unique numbers $a,b,c>0$ such that
$abi(P_\alpha,P_\beta)=bci(P_\beta,P_\gamma)=aci(P_\gamma,P_\alpha)=1$.
By construction, we have
\[L_a(P_\alpha,P_\beta,k_1)=L_b(P_\beta,P_\alpha,k_1),\,
L_b(P_\beta,P_\gamma,k_1)=L_c(P_\gamma,P_\beta,k_1)\] and
$L_c(P_\gamma,P_\alpha,k_1)=L_a(P_\alpha,P_\gamma,k_1)$.
Choose a point $\delta\in L_a(P_\alpha,P_\beta,k_1)$ such that
for every $\zeta$ in ${\cal M\cal C}(S)$ we have
\[i(\delta,\zeta)\leq k_1\max\{ai(P_\alpha,\zeta),
i(\zeta,P_\beta)/ai(P_\alpha,P_\beta)\};\] such a point exists by Lemma
3.3. Applying this inequality to $\zeta=P_\gamma$ yields
$ci(\delta,P_\gamma)\leq k_1$. For $\zeta=P_\beta$ we obtain
\[i(\delta,P_\beta)/ci(P_\gamma,P_\beta)\leq
ai(P_\alpha,P_\beta)/ci(P_\gamma,P_\beta)=1,\]
and for $\zeta=P_\alpha$ we obtain
\[i(\delta,P_\alpha)/ci(P_\gamma,P_\alpha)
\leq 1/aci(P_\gamma,P_\alpha)=1.\]
Therefore we have $\delta\in
L_c(P_\gamma,P_\beta,k_1)\cap L_a(P_\alpha,P_\beta,k_1)
\cap L_c(P_\gamma,P_\alpha,k_1)$.
Together with Lemma 3.2 and our above remark we conclude that
there is a universal constant $D_2>0$ such that the distance
between
$\phi(\alpha,\beta,\gamma)=\delta$ and the
image under $\Phi$ of any $\alpha\to\beta$-splitting sequence,
any $\alpha\to\gamma$-splitting sequence and any
$\gamma\to\beta$-splitting
sequence is bounded from above by $D_2$.

We use the map $\phi$ to derive the $\delta$-thin
triangle condition for triangles whose sides
are images under the map $\Phi$ of splitting sequences
in ${\cal T}T$.

\begin{lemma}
There is a number $D_3>0$ with the following
property.
Let $\alpha,\beta,\gamma\in {\cal C}(S)$
and let
$a,b,c$ be the image under $\Phi$ of a $\beta\to\gamma$,
$\gamma\to \alpha$, $\alpha\to \beta$-splitting sequence.
Then the $D_3$-neighborhood of $a\cup b$ contains
$c$.
\end{lemma}

\begin{proof} Let $\alpha,\beta,\gamma\in {\cal C}(S)$ and assume that
$d(\beta,\gamma)\leq p$ for some $p>0$.
Let $(\tau_i)_{0\leq i\leq m}$ be an $\alpha\to \beta$-splitting
sequence and let $(\eta_j)_{0\leq j\leq \ell}$
be an $\alpha\to\gamma$-splitting
sequence; if $D_1>0$ is as above
then the Hausdorff distance
between $\{\Phi(\tau_i)\mid 0\leq i\leq m\}$ and
$\{\Phi(\eta_j)\mid 0\leq j\leq \ell\}$ is at most $2pD_1$.
Namely, we observed that the Hausdorff distance
between the image under $\Phi$ of any two
$\alpha\to\beta$-splitting sequences is bounded
from above by $D_1$. Moreover,
if $d(\beta,\gamma)=1$ then $\beta\cup\gamma\in
{\cal M\cal C}(S)$ and hence
there is an $\alpha\to\beta$-splitting
sequence which also is an $\alpha\to\gamma$-splitting sequence.
Thus the statement of the corollary holds for
$p=1$, and the
general case follows from a successive application of this fact for
the points on a geodesic 
in ${\cal C}(S)$ connecting $\beta$ to $\gamma$.

Now let $\alpha,\beta,\gamma\in {\cal C}(S)$
be arbitrary points whose
pairwise distance is at least 3. 
Let again $(\tau_i)_{0\leq i\leq m}$ be
an $\alpha\to\beta$-splitting sequence.
By the definition of $\phi$ and the choice of the
constant $D_2>0$ above
there is
some $i_0\leq m$ such that the distance
betweeen $\Phi(\tau_{i_0})$ and $\phi(\alpha,\beta,\gamma)$ is
at most $D_2$. Let $(\eta_j)_{0\leq j\leq \ell}$ be any
$\alpha\to\phi(\alpha,\beta,\gamma)$-splitting sequence.
By our above consideration, the Hausdorff distance
between $\{\Phi(\tau_i)\mid 0\leq i\leq i_0\}$ and
$\{\Phi(\eta_j)\mid 0\leq j\leq \ell\}$ is at most $2D_1D_2$.
Similarly, let $(\zeta_j)_{0\leq j\leq n}$ be
an $\alpha\to\gamma$-splitting sequence.
Then there is some $j_0>0$ such that $d(\Phi(\zeta_{j_0}),
\phi(\alpha,\beta,\gamma))\leq D_2$. By our above argument,
the Hausdorff distance between the sets
$\{\Phi(\tau_i)\mid 0\leq i\leq i_0\}$ and
$\{\Phi(\zeta_j)\mid 0\leq j\leq j_0\}$ is at most
$4D_1D_2$.

As a consequence, there are numbers
$a(\alpha,\beta)>0,a(\alpha,\gamma)>0$ such that
\[\phi(\alpha,\beta,\gamma)\in 
L_{a(\alpha,\beta)}(\beta,\alpha,k_1)\cap
L_{a(\alpha,\gamma)}(\gamma,\alpha,k_1)\]
and that the Hausdorff distance between
$\cup_{a\geq a(\alpha,\beta)}L_a(\beta,\alpha,k_1)$
and\\ $\cup_{a\geq a(\alpha,\gamma)}L_a(\gamma,\alpha,k_1)$
is at most $6D_1D_2$. The same argument, applied
to a $\beta\to \alpha$-splitting sequence
and a $\beta\to\gamma$-splitting sequence,
shows that $\cup_{a\leq a(\alpha,\beta)}L_a(\alpha,\beta,k_1)$
is contained in the $6D_1D_2$-neighborhood of
$\cup_aL_a(\beta,\gamma,k_1)$. Then
$\{\Phi(\tau_i)\mid 0\leq i\leq m\}$ is contained
in the $12D_1D_2$-neighborhood of the union of the image
under $\Phi$ of a $\gamma\to\alpha$-splitting sequence
and a $\beta\to\gamma$-splitting sequence. This shows the lemma.
\end{proof}

Hyperbolicity of the curve graph
now follows from Lemma 3.4 and the
following criterion.

\begin{proposition}
Let $(X,d)$ be a geodesic metric space.
Assume that there is a number $D>0$ and for
every pair of points $x,y\in X$ there is an
arc $\eta(x,y):[0,1]\to X$ connecting $\eta(x,y)(0)=x$
to $\eta(x,y)(1)=y$ so that the following conditions
are satisfied.
\begin{enumerate}
\item If $d(x,y)\leq 1$ then the diameter of 
$\eta(x,y)[0,1]$ is at most $D$.
\item For $x,y\in X$ and
$0\leq s\leq t\leq 1$, the Hausdorff
distance between $\eta(x,y)[s,t]$ and 
$\eta(\eta(x,y)(s),\eta(x,y)(t))[0,1]$ is at most $D$.
\item For any $x,y,z\in X$ the set $\eta(x,y)[0,1]$
is contained in the $D$-neighborhood of 
$\eta(x,z)[0,1]\cup \eta(z,y)[0,1].$
\end{enumerate}
Then $(X,d)$ is $\delta$-hyperbolic for a number $\delta>0$
only depending on $D$.
\end{proposition}

\begin{proof} Let $(X,d)$ be a geodesic
metric space. Assume that there is a number
$D>0$ and there is a family of paths $\eta(x,y):[0,1]\to X$,
one for every pair of points $x,y\in X$, which
satisfy the hypotheses in the statement of the proposition.
To show hyperbolicity for $X$ it is then enough to 
show the existence
of a constant $\kappa >0$
such that for all $x,y\in X$ and
every geodesic 
$\nu:[0,\ell]\to X$
connecting $x$ to $y$, the Hausdorff-distance
between $\nu[0,\ell]$ and
$\eta(x,y)[0,1]$ is at most $\kappa$.
Namely, if this is the case then
for every
geodesic triangle with sides $a,b,c$ the side
$a$ is contained in the $3\kappa+D$-neighborhood of $b\cup c$.

To show the existence of such a constant $\kappa >0$,
let $x,y\in X$ and 
let $c:[0,2^k]\to X$ be \emph{any}
path of length $\ell(c)=2^k$ parametrized by arc length
connecting $x$ to $y$. 
Write $\eta_1=\eta(c(0),c(2^{k-1}))$ and write $\eta_2=
\eta(c(2^{k-1}),c(2^k)).$
By our assumption, the
$D$-neighborhood of $\eta_1\cup \eta_2$ contains
$\eta(c(0),c(2^k))$.
Repeat this construction with the points
$c(2^{k-2}),c(3\cdot 2^{k-2})$ and the arcs
$\eta_1,\eta_2$. Inductively we conclude that the path
$\eta(c(0),c(2^k))$ is contained in the
$(\log_2\ell(c))D$-neighborhood of a path 
$\tilde c:[0,2^k]\to {\cal C}(S)$ whose restriction
to each interval $[m-1,m]$ $(m\leq 2^k)$ equals 
up to parametrization the arc
$\eta(c(m-1),c(m))$. Since $d(c(m-1),c(m))\leq 1$, by
assumption the diameter of each of the sets $\eta(c(m-1),c(m))[0,1]$
is bounded from above by $D$ and therefore
the arc $\eta(c(0),c(2^k))$
is contained in the $(\log_2\ell(c))D+D$-neighborhood 
of $c[0,2^k]$.

Now let $c:[0,k]\to X$ be a geodesic connecting
$c(0)=x$ to $c(k)=y$ which
is parametrized by arc length. Let $t>0$ be such that
$\eta(x,y)(t)$ has maximal distance to $c[0,k]$,
say that this distance equals $\chi$. 
Choose some $s>0$ such that
$d(c(s),\eta(x,y)(t))=\chi$ and let $t_1<t<t_2$ be such that
$d(\eta(x,y)(t),\eta(x,y)(t_u))= 2\chi$ $(u=1,2)$.
In the case that there is no
$t_1\in [0,t)$ (or $t_2\in (t,1]$)
with $d(\eta(x,y)(t),\eta(x,y)(t_1))\geq 2\chi$
(or $d(\eta(x,y)(t),\eta(x,y)(t_2))\geq 2\chi$) we
choose $t_1=0$ (or $t_2=1$).
By our choice of $\chi$, there 
are numbers $s_u\in
[0,k]$ such that $d(c(s_u),\eta(x,y)(t_u))\leq \chi
\,(u=1,2)$.
Then the distance
between $c(s_1)$ and $c(s_2)$ is at most $6\chi$. Compose the
subarc $c[s_1,s_2]$ of $c$
with a geodesic
connecting $\eta(x,y)(t_1)$ to $c(s_1)$ and a geodesic
connecting $c(s_2)$ to $\eta(x,y)(t_2)$. We obtain a curve $\nu$
of length at most $8\chi$. By our above observation, the
$(\log_2(8\chi))D+D$-neighborhood of this curve contains the
arc $\eta(\eta(x,y)(t_1),\eta(x,y)(t_2))$.
However, the Hausdorff distance between
$\eta(x,y)[t_1,t_2]$ and $\eta(\eta(x,y)(t_1),
\eta(x,t)(t_2))$ is at most $D$ and therefore the
$(\log_2(8\chi))D+2D$-neighborhood of the arc $\nu$ contains
$\eta(x,y)[t_1,t_2]$. 
But the distance between
$\eta(t)$ and our curve $\nu$ equals $\chi$ by construction
and hence we have 
$\chi\leq (\log_2(8\chi))D+2D$. In other words, $\chi$
is bounded from above
by a universal constant $\kappa_1 >0$, and
$\eta(x,y)$
is contained in the $\kappa_1$-neighborhood of
the geodesic $c$.

A similar argument also shows that the 
$3\kappa_1$-neighborhood of 
$\eta(x,y)$ contains $c[0,k]$. Namely, 
by the above consideration, for every $t\leq 1$ 
the set $A(t)=\{s\in [0,k]\mid 
d(c(s),\eta(x,y)(t))\leq \kappa_1\}$
is a non-empty closed subset of $[0,k]$.
The diameter of the sets $A(t)$ is bounded
from above by $2\kappa_1$. 
Assume to the contrary that 
$c[0,k]$ is not contained in the 
$3\kappa_1$-neighborhood of $\eta(x,y)$.
Then there is a subinterval
$[a_1,a_2]\subset [0,k]$ of length $a_2-a_1\geq 4\kappa_1$ 
such that $d(c(s),\eta(x,y)[0,1])>\kappa_1$ for
every $s\in (a_1,a_2)$. Since $d(c(s),c(t))=\vert s-t\vert$
for all $s,t$ we conclude that for every
$t\in [0,1]$ the set $A(t)$ either is entirely
contained in $[0,a_1]$ or it is entirely contained
in $[a_2,k]$. Define $C_1=\{t\in [0,1]\mid
A(t)\subset [0,a_1]\}$ and $C_2=\{t\in [0,1]\mid
A(t)\subset [a_2,1]\}$. Then the sets $C_1,C_2$ are disjoint and
their union equals $[0,1]$; moreover, we
have $0\in C_1$ and $1\in C_2$. On the other
hand, the sets $C_i$ are closed. Namely,
let $(t_i)\subset C_1$ be a sequence converging
to some $t\in [0,1]$. 
Let $s_i\in A(t_i)$ and assume after passing
to a subsequence that $s_i\to s\in [0,a_1]$.
Now $\kappa_1\geq d(c(s_i),\eta(x,y)(t_i))\to
d(c(s),\eta(x,y)(t))$ and therefore $s\in A(t)$
and hence $t\in C_1$. However, $[0,1]$ is connected
and hence we arrive at a contradiction. 
In other words, the geodesic $c$ is contained
in the $3\kappa_1$-neighborhood of $\eta(x,y)$.
This completes the proof of the proposition.
\end{proof}

As an immediate corollary, we obtain.

\begin{theorem}[\cite{MM1,B}]
The curve graph is hyperbolic.
\end{theorem}

\begin{proof} Let $\Phi:{\cal T}T\to {\cal C}(S)$
be as before. For $\alpha,\beta\in {\cal C}(S)$ choose
an $\alpha\to\beta$-splitting sequence
$(\tau_i)_{0\leq i\leq m}$. 
Define an arc $\eta(\alpha,\beta):[0,1]\to
{\cal C}(S)$ by requiring that for $1\leq i\leq m$ 
the restriction of 
$\eta(\alpha,\beta)$ to the interval
$[\frac{i}{m+2},\frac{i+1}{m+2}]$ is a geodesic
connecting $\Phi(\tau_{i-1})$ to 
$\Phi(\tau_i)$ and that the restriction of $\eta(\alpha,\beta)$
to $[0,\frac{1}{m+2}]$ (or $[\frac{m-1}{m-2},1]$)
is a geodesic connecting $\alpha$ to 
$\Phi(\tau_0)$ (or $\Phi(\tau_m)$ to $\beta$).

We claim that this family of arcs satisfy the assumptions
in Proposition 3.5. Namely, we observed before
that for all $\alpha,\beta\in {\cal C}(S)$,
the Hausdorff distance between the image
under $\Phi$ of \emph{any} two $\alpha\to\beta$-splitting
sequences is bounded from above by a universal constant.
Now if $(\tau_i)_{0\leq i\leq m}$ is any splitting
sequence, then for all $0\leq k\leq \ell\leq m$ the
sequence $(\tau_i)_{k\leq i\leq \ell}$ is
a $\Phi(\tau_k)\to\Phi(\tau_\ell)$-splitting sequence
and hence our curve system satisfies the
second condition in Proposition 3.5.

Moreover, curves $\alpha,\beta\in {\cal C}(S)$ 
with $d(\alpha,\beta)=1$ can be realized disjointly,
and $\alpha\cup \beta$ is a multi-curve.
For such a pair of curves we can
choose a \emph{constant} $\alpha\to\beta$-splitting
sequence; hence our curve system also
satisfies the first condition stated in Proposition 3.5.
Finally, the third condition was shown to hold in 
Lemma 3.4. 

Now hyperbolicity of the curve graph
follows from Proposition 3.5.
\end{proof}

A curve $c:[0,m]\to {\cal C}(S)$ is called an \emph{unparametrized
$p$-quasi-geodesic}\index{unparametrized
quasi-geodesic} for some $p>1$ if there is a homeomorphism
$\rho:[0,u]\to [0,m]$ for some $u>0$ such that
\[d(c(\rho(s)),c(\rho(t)))/p-p\leq \vert s-t\vert \leq
pd(c(\rho(s)),c(\rho(t)))+p\] for all $s,t\in [0,u]$.
We define a map $c:\{0,\dots,m\}\to {\cal C}(S)$ to be
an unparametrized $q$-quasi-geodesic
if this is the case for the curve $\tilde c$ whose
restriction to each interval $[i,i+1)$ coincides with
$c(i)$.  
The following observation is immediate from Proposition 3.5
and its proof.

\begin{corollary}[\cite{MM3, H2}]
There is a number $p>0$
such that the image under $\Phi$ of an arbitrary splitting
sequence is an unparametrized $p$-quasi-geodesic.
\end{corollary}

\begin{proof} By Proposition 3.5, 
Theorem 3.6 and their proofs, there is
a universal number $D>0$ with the property that
for every splitting sequence $(\tau_i)_{0\leq i\leq m}$ and
every geodesic $c:[0,s]\to {\cal C}(S)$ connecting
$c(0)=\Phi(\tau_0)$ to $c(s)=\Phi(\tau_m)$, the Hausdorff distance
between the sets $\{\Phi(\tau_i)\mid 0\leq i\leq m\}$ and
$c[0,s]$ is at most $D$. From this the
corollary is immediate.
\end{proof}

\section{Geometry of Teichm\"uller space}

In this section, we relate the geometry of the curve
graph to the geometry of Teichm\"uller space
equipped with the Teichm\"uller metric. For this
we first define a map $\Psi:{\cal T}_{g,m}\to {\cal C}(S)$
as follows.
By a well-known result of Bers (see \cite{Bu}) there is a number
$\chi>0$ only depending on the topological type of $S$
such that
for every complete hyperbolic metric 
on $S$ of finite volume
there is a pants decomposition $P$ for $S$ which consists
of simple closed geodesics of length at most $\chi$.
Since the distance between any two points $\alpha,\beta
\in {\cal C}(S)$ is bounded from above by
$i(\alpha,\beta)+1$, the collar lemma for hyperbolic
surfaces (see \cite{Bu}) implies that 
the diameter in ${\cal C}(S)$ of the set
of simple closed curves whose length with respect to
the fixed metric is at most $\chi$ is bounded
from above by a universal constant $D>0$. Define
$\Psi:{\cal T}_{g,m} \to
{\cal C}(S)$ by assigning to a finite volume hyperbolic
metric $h$ on $S$ a simple closed curve $\Psi(h)$ whose
$h$-length is at most $\chi$. Then
for any two maps
$\Psi,\Psi^\prime$ with this property and every $h\in {\cal T}_{g,m}$
the distance in ${\cal C}(S)$ between $\Psi(h)$ and
$\Psi^\prime(h)$ is at most $D$. Moreover, the map $\Psi$
is coarsely equivariant with respect to the
action of the mapping class group ${\cal M}_{g,m}$
on ${\cal T}_{g,m}$ and ${\cal C}(S)$: For every
$h\in {\cal T}_{g,m}$ and every $\phi\in {\cal M}_{g,m}$ we have
$d(\Psi(\phi(h)),\phi(\Psi(h)))\leq D$.

The following result is due to Masur and Minsky 
(Theorem 2.6 and Theorem 2.3 of \cite{MM1}).
For its formulation, let $d_T$ be the distance function
on ${\cal T}_{g,m}$ induced by the Teichm\"uller metric.

\begin{theorem}[\cite{MM1}]
\begin{enumerate}
\item There is a number $a>0$ such that
$d(\Psi h,\Psi h^\prime)\leq ad_T(h,h^\prime)+a$ for
all $h,h^\prime\in {\cal T}_{g,m}$.
\item
There is a number $\tilde p>0$
with the following property.
Let
$\gamma:(-\infty,\infty)\to {\cal T}_{g,m}$ be any
Teichm\"uller geodesic; then the
assignment $t\to \Psi(\gamma(t))$ is 
an unparametrized $\tilde p$-quasi-geodesic
in ${\cal C}(S)$.
\end{enumerate}
\end{theorem}

\begin{proof} 
Let $\gamma:(-\infty,\infty)\to {\cal T}_{g,m}$
be any Teichm\"uller geodesic parametrized
by arc length. Then the cotangent of $\gamma$
at $t=0$ is a quadratic differential $q$
of area one defined by a pair  
$(\lambda,\mu)\in {\cal M\cal L}\times {\cal M\cal L}$ of 
measured geodesic laminations which jointly fill
up $S$. The cotangent of $\gamma$ at $t$ is given 
by the quadratic differential $q(t)$
defined by the pair $(e^t\lambda,e^{-t}\mu)$.
For $k_1>0$ as in Lemma 3.3 and $t\in \mathbb{R}$ 
let $\zeta(t)\in {\cal C}(S)$ be a curve whose
$q(t)$-length is at most $2k_1$. For every $\beta\in [0,1]$
the $q(t+\beta)$-length
of $\zeta(t)$ is bounded from above by $2ek_1$ and
therefore by Lemma 3.3, for every $t$ the distance
in ${\cal C}(S)$ between $\zeta(t)$ and $\zeta(t+\beta)$
is bounded from above by a universal constant $k_2>0$.
In particular, the assignment $t\to \zeta(t)$
satisfies $d(\zeta(s),\zeta(t))\leq k_2\vert s-t\vert 
+k_2$. Hence for the proof of our lemma, we only have to show
that there is a constant $k_3>0$ such that for
every $h\in {\cal T}_{g,m}$ and every holomorphic
quadratic differential $q$ of area one for $h$,
the distance between $\Psi(h)$ and a curve on $S$
whose $q$-length is bounded from above by $2k_1$ is 
uniformly bounded.

Thus let $h$ be a complete hyperbolic metric of finite
volume and let $q$ be a holomorphic 
quadratic differential for
$h$ of area one. By the collar lemma of hyperbolic geometry,
a simple closed geodesic $c$ for $h$ whose
length is bounded from above by $\chi$ is
the core curve of an embedded annulus $A$ whose
\emph{modulus}\index{modulus} is bounded from
below by a universal constant $\epsilon >0$; we
refer to \cite{S} for a definition of the modulus
of an annulus and its properties.
Then the \emph{extremal length} of the core curve
of $A$ is bounded from above by a universal
constant $m >0$. Now the area of 
$q$ equals one and therefore
the $q$-length of the core curve $c$ does not
exceed $\sqrt{m}$ by the definition of
extremal length (see e.g. \cite{Mi}).
In other words, the $q$-length
of the curve $\Psi(h)$ is uniformly bounded
which together with Lemma 3.3 implies our 
claim. The theorem follows.
\end{proof}

There are also Teichm\"uller geodesics
in Teichm\"uller space which are mapped by
$\Psi$ to \emph{parametrized} quasi-geodesics in 
${\cal C}(S)$. For their characterization, denote for
$\epsilon >0$ by ${\cal T}_{g,m}^\epsilon$ the subset
of Teichm\"uller space consisting of all marked hyperbolic
metrics for which the length of the shortest closed
geodesic is at least $\epsilon$. The set ${\cal T}_{g,m}^\epsilon$
is invariant under the action of the mapping class group
and projects to a \emph{compact} subset of moduli
space. 
Moreover, every compact subset of 
moduli space is contained in the projection of 
${\cal T}_{g,m}^\epsilon$ for some $\epsilon >0$.

\emph{Cobounded} Teichm\"uller geodesics. i.e. Teichm\"uller
geodesic 
which project to a compact subset of moduli space, 
relate the geometry of Teichm\"uller space
to the geometry of the curve graph.
We have.

\begin{proposition}[\cite{H05}]
The image under $\Psi$ of a Teichm\"uller geodesic
$\gamma:\mathbb{R}\to {\cal T}_{g,m}$ is a
parametrized quasi-geodesic in ${\cal C}(S)$ if and 
only if there is some $\epsilon >0$ such that
$\gamma(\mathbb{R})\subset {\cal T}_{g,m}^\epsilon$.
\end{proposition}

Minsky \cite{Mi96} discovered earlier that 
the Teichm\"uller
metric near a cobounded geodesic line 
has properties similar to properties of a hyperbolic
geodesic metric space. Namely, for a Teichm\"uller geodesic
$\gamma:\mathbb{R}\to {\cal T}_{g,m}^\epsilon$ the 
map which associates to a point $h\in {\cal T}_{g,m}$
a point on $\gamma(\mathbb{R})$ which
minimizes the Teichm\"uller distance is coarsely
Lipschitz and contracts distances
in a way which is similar to the contraction
property of the
closest point projection 
from a $\delta$-hyperbolic geodesic metric space to
any of its bi-infinite geodesics.

A hyperbolic geodesic metric space $X$ admits a \emph{Gromov
boundary}\index{Gromov boundary} 
which is defined as follows. Fix a point $p\in X$ and
for two points $x,y\in X$ define the \emph{Gromov product}
$(x,y)_p=\frac{1}{2}(d(x,p)+d(y,p)-d(x,y))$. Call a sequence
$(x_i)\subset X$ \emph{admissible} if $(x_i,x_j)_p\to \infty$
$(i,j\to \infty)$. We define two admissible sequences
$(x_i),(y_i)\subset X$ to be \emph{equivalent} if $(x_i,y_i)_p\to
\infty$. Since $X$ is hyperbolic, this defines indeed an
equivalence relation (see \cite{BH}). The Gromov boundary $\partial X$
of $X$ is the set of equivalence classes of admissible
sequences $(x_i)\subset X$. It carries a natural Hausdorff
topology. For the curve graph,
the Gromov boundary was determined by Klarreich \cite{K} (see
also \cite{H2}).

For the formulation of Klarreich's result,
we say that a minimal geodesic
lamination $\lambda$ \emph{fills up $S$}
if every simple closed geodesic on $S$ intersects $\lambda$
transversely, i.e. if every complementary component of $\lambda$
is an ideal polygon or a once punctured ideal polygon with
geodesic boundary \cite{CEG}.
For any minimal geodesic lamination $\lambda$ which fills up $S$,
the number of geodesic laminations $\mu$ which contain $\lambda$
as a sublamination is bounded by a universal constant only
depending on the topological type of the surface $S$. Namely, each
such lamination $\mu$ can be obtained from $\lambda$ by
successively subdividing complementary components $P$ of $\lambda$
which are different from an ideal triangle or a once punctured
monogon by adding a simple geodesic line which either
connects two non-adjacent cusps or goes around a puncture.
Note that every leaf of $\mu$ which is not
contained in $\lambda$ is necessarily isolated in $\mu$.

Recall that the space ${\cal L}$ 
of geodesic laminations on $S$ equipped with the
restriction of the 
Hausdorff topology for compact subsets of $S$
is compact and metrizable. It contains the 
set ${\cal B}$ of minimal geodesic laminations
which fill up $S$ as a subset which is neither closed
nor dense. We define on ${\cal B}$ a new topology
which is coarser than the restriction of the Hausdorff topology 
as follows.
Say that a sequence $(\lambda_i)\subset {\cal L}$
\emph{converges in the coarse Hausdorff topology} to a
minimal geodesic lamination $\mu$ which fills up $S$ if
every accumulation point of $(\lambda_i)$ with respect
to the Hausdorff topology contains $\mu$ as a sublamination.
Define a subset $A$ of 
${\cal B}$ to be closed if and only if
for every sequence
$(\lambda_i)\subset A$ which converges in the coarse
Hausdorff topology to a lamination $\lambda\in {\cal B}$
we have $\lambda\in A$. 
We call the resulting topology on
${\cal B}$ the 
\emph{coarse Hausdorff topology}\index{coarse Hausdorff topology}.
The space ${\cal B}$ is
not locally compact.
Using this terminology, Klarreich's
result \cite{K} can be formulated as follows.

\begin{theorem}[\cite{K, H2}]
\begin{enumerate}
\item There is a natural homeomorphism $\Lambda$ of ${\cal B}$\\
equipped with the coarse Hausdorff topology onto the Gromov
boundary $\partial {\cal C}(S)$ of the complex of curves
${\cal C}(S)$ for $S$.
\item For $\mu\in {\cal B}$
a sequence $(c_i)\subset {\cal C}(S)$ is admissible and
defines the point $\Lambda(\mu)\in \partial {\cal C}(S)$ if and
only if $(c_i)$ converges in the coarse Hausdorff
topology to $\mu$.
\end{enumerate}
\end{theorem}

Recall that every Teichm\"uller geodesic in ${\cal T}_{g,m}$
is uniquely determined by a pair of projective
measured laminations which jointly fill up $S$. 
The following corollary is immediate from Theorem 4.1 and 
Theorem 4.3 with $\tilde p>0$ as in Theorem 4.1.

\begin{corollary}
Let $\lambda,\mu\in {\cal P\cal M\cal L}$ be such that
their supports are minimal and fill up $S$. Then the
image under $\Psi$ of the
unique Teichm\"uller geodesic in ${\cal T}_{g,m}$ 
determined by $\lambda$ and $\mu$ is a biinfinite
unparametrized $\tilde p$-quasigeodesic in ${\cal C}(S)$, 
and every biinfinite unparametrized $\tilde p$-quasi-geodesic
in ${\cal C}(S)$ is contained in a uniformly bounded
neighborhood of a curve of this form. 
\end{corollary}

Our above discussion also gives information on images 
under $\Psi$ of a convergent sequence of 
geodesic lines in Teichm\"uller space.
Namely, if $(\gamma_i)$ is such a sequence of Teichm\"uller
geodesic lines converging to
a Teichm\"uller geodesic which is determined by a pair of
projective measured geodesic laminations $(\alpha,\beta)$ so that
the support of $\alpha$ does not fill up $S$ then
there is a curve $\zeta\in {\cal C}(S)$, a 
number $m>0$ and a sequence
$j(i)\to \infty$ such that $d(\Psi(\gamma_i[0,j(i)]),\zeta)\leq m$.

On the other hand, the image under $\Psi$ of ``most''
Teichm\"uller geodesics
are unparametrized quasi-geodesics
of infinite diameter which are not \emph{parametrized}
quasi-geodesics. For example, let $\lambda\in {\cal P\cal M\cal L}$
be a projective measured geodesic lamination whose
support $\lambda_0$ is minimal and fills up $S$ 
but is not \emph{uniquely
ergodic}. This means that the dimension of the
space of transverse measures supported in $\lambda_0$
is at least 2. Let $\gamma:[0,\infty)\to
{\cal T}_{g,m}$ be a Teichm\"uller geodesic
ray determined by a quadratic differential whose horizontal
foliation corresponds to $\lambda$. By a result of Masur
\cite{Ma82a}, the projection of $\gamma$ to moduli
space eventually leaves every compact set.
On the other hand,
since $\lambda_0$ is minimal and fills up $S$ the points
$\gamma(t)$ converge as $t\to \infty$ to $\lambda$ viewed
as a point in the \emph{Thurston boundary} of 
the Thurston compactification of Teichm\"uller space
\cite{Ma82b} (compare also \cite{FLP} 
for the construction of the Thurston compactification).
By the definition of the map $\Psi$, the projective
measured geodesic laminations defined by the curves
$\Psi\gamma(t)$ converge as $t\to\infty$ to $\lambda$ and
therefore the curves $\Psi\gamma(t)$ converge in the
coarse Hausdorff topology to $\lambda_0$. By Theorem 4.3,
this implies that the diameter of $\Psi\gamma[0,\infty)$ 
is infinite.

\end{document}